\documentclass[9pt,final,twocolumn]{IEEEtran}

\newtheorem{prob}{Problem}
\newtheorem{prop}{Proposition}
\newtheorem{mydef}{Definition}
\newtheorem{lem}{Lemma}
\newtheorem{cor}{Corollary}

\newtheorem{theorem}{Theorem}
\newtheorem{rem}{Remark}

\usepackage{cite}

\ifCLASSINFOpdf
\usepackage[pdftex]{graphicx}
\graphicspath{{../pdf/}{../jpeg/}}
\DeclareGraphicsExtensions{.pdf,.jpeg,.png}
\else
\usepackage[dvipdf]{graphicx}
\graphicspath{{../eps/}}
\DeclareGraphicsExtensions{.eps}
\fi

\usepackage[cmex10]{amsmath}
\usepackage{amssymb}
\usepackage{array}
\usepackage{mdwmath}
\usepackage{mdwtab}
\usepackage{eqparbox}
\usepackage[tight,footnotesize]{subfigure}

\begin{document}

\title{Minimum-Time Transitions between Thermal Equilibrium States of the Quantum Parametric Oscillator}

\author{Dionisis~Stefanatos,~\IEEEmembership{Member,~IEEE}
\thanks{D. Stefanatos is with the Division of Physical Sciences and Applications, Hellenic Army Academy, Vari, Athens 16673, Greece, e-mail: dionisis@post.harvard.edu.}
}

\markboth{}
{Stefanatos: Minimum-Time Thermal Equilibrium Transitions}

\maketitle

\begin{abstract}

In this article, we use geometric optimal control to completely solve the problem of minimum-time transitions between thermal equilibrium states of the quantum parametric oscillator, which finds applications in various physical contexts. We discover a new kind of optimal solutions, absent from all the previous treatments of the problem.

\end{abstract}


\IEEEpeerreviewmaketitle

\section{Introduction}

\IEEEPARstart{C}{onstantin} Carath\'{e}odory, the famous mathematician with seminal contributions in the calculus of variations \cite{Caratheodory35} which paved the way to optimal control theory \cite{Pesch13}, pioneered the axiomatic formulation of thermodynamics along a purely geometric approach \cite{Caratheodory09}, at the dawn of the 20th century. Based on these foundations, the geometry of thermodynamics \cite{Weinhold09} was developed many decades later. In this context, thermal equilibrium states are represented as points on a manifold and tools from differential geometry are used to quantify the distance between them and to express the laws of thermodynamics. Closely related to this approach is the subject of finite-time thermodynamics \cite{Andresen84,Berry99}, which aims to optimize the performance of a thermodynamic system under restrictions on the available time, for example to maximize the extracted power. Optimal control theory \cite{Pontryagin} is the mathematical tool used to tackle this kind of problems, thus the connection between thermodynamics and control is not just restricted to the emblematic figure of Carath\'{e}odory but is actually deeper.

The field of finite-time thermodynamics has been recently revitalized in the context of quantum systems, where the design and optimization of nanoscale heat engines provides the major motivation \cite{Feldman96,Salamon09,Abah12,Stefanatos14PRE,Bonanca14,Azimi14,Campo14,Rotskoff15,Xiao15,Zulkowski15}. One paradigmatic example which stands out from the rest is the quantum parametric oscillator \cite{Salamon09}, a quantum harmonic oscillator whose angular frequency can be altered with time and serves as the control parameter. The problem related to this system is, starting from a thermal equilibrium state and changing the frequency from some initial to a lower final value, to find the maximum work that can be extracted and the minimum necessary time. The authors of \cite{Salamon09} consider the realistic case where the frequency of the oscillator can only take real values, corresponding to nonnegative stiffness, which is restricted between a lower and an upper bound. Under these assumptions, they show that the maximum work is obtained when the final state of the system is also an equilibrium state, and provide an analytical estimate of the minimum time, which depends on the frequency bounds. The minimum-time solution is achieved when the frequency changes in a bang-bang manner between its boundary values following a three-jumps strategy, including two jumps at the initial and final times and only one \emph{intermediate} switching.

The importance of this problem stems from its applicability in various frameworks. First, as a minimum-time problem between two thermal states, it can find application is solid-state chemistry, in the transition from graphite to diamond, as pointed out by the authors in \cite{Salamon09}. Second, if the frequency at the final state is lower than the frequency at the initial state, the temperature corresponding to the final thermal state is lower than the initial temperature and this corresponds to cooling the system. During the process of frequency change the system is isolated from its environment, so its entropy remains constant and the process is adiabatic in the thermodynamic sense. Thus, the solution of the problem provides also the minimum time for an adiabatic stroke in a quantum heat engine, operating between a hot and a cold reservoir. But the described process is \emph{effectively} adiabatic also in the quantum sense, since the populations of the energy levels at the final state are the same with those of the corresponding energy levels at the initial state. Recall from elementary quantum physics that for the process to be adiabatic it would be necessary the populations to remain the same during the whole procedure, and this would require a slow enough variation of the frequency from the initial to the final value. In contrast, during the minimum-time process the frequency changes abruptly and only at the final time is recovered the initial population distribution, thus it is characterized as \emph{effectively} adiabatic. In conclusion, the minimum-time solution provides also the fastest effectively adiabatic cooling for a particle trapped in the parametric harmonic oscillator, a procedure with many interesting applications in physics \cite{Chen10}.

It was in this last context of adiabatic-like cooling where ``shortcuts to adiabaticity" were introduced \cite{Chen10}, protocols where the time-profile of the frequency is obtained by appropriate interpolation between the initial and final values. The authors of \cite{Chen10} permitted the possibility that the harmonic potential can become expulsive for some time intervals, which corresponds to imaginary frequency and negative stiffness of the oscillator (the actual control in this case), and did not apply any control bounds. Under these broader assumptions, compared to those in the original formulation \cite{Salamon09} where only real bounded frequencies were allowed, they concluded that the desired transition between the initial and final thermal states can take place in arbitrarily short times. By imposing bounds on the stiffness of the oscillator, which are always present because of experimental limitations, we were able to show that even in the case where a expulsive potential is allowed, there is a minimum necessary time for the transition between the thermal states \cite{Stefanatos10PRA}. We formulated the corresponding optimal control problem and proved that the optimal solution has the bang-bang form, as in the original case \cite{Salamon09}. In our subsequent work \cite{Stefanatos11} we completely solved the problem where negative stiffness values are allowed, and obtained a type of solution with an even number of intermediate switchings, which was absent from the original work \cite{Salamon09}. Recently, it was numerically confirmed \cite{Boldt16} by some authors of \cite{Salamon09} that our solution is encountered even in the more restrictive case, for appropriate values of the parameters (control bounds, ratio of initial and final frequency).

This last work gave as the motivation to rigorously study the optimal control problem for the restrictive case with nonnegative stiffness. As the authors in \cite{Hoffmann13} point out, ``more restrictive controls can lead to more interesting answers that reveal more of the physics of the problem". Note that optimal control theory \cite{Pontryagin} is one of the basic methods of quantum control \cite{Altafini12} and has been successfully applied to obtain minimum-time solutions for several quantum systems \cite{Khaneja01,Boscain06,Bonnard09,Stefanatos11,Bonnard12,Stefanatos13,Stefanatos14TAC,Albertini14,Zhang15}, in an attempt to reduce the undesirable interactions with the environment which lead to dissipation and decoherence. In the present article, we use geometric optimal control \cite{Heinz12} and completely solve the problem of minimum-time transitions between thermal states of the quantum parametric oscillator, for the restrictive case where the stiffness can only take nonnegative values (real frequency) \cite{Salamon09}. We recover the one intermediate switching solution, presented in \cite{Salamon09}, and our solution with even number of intermediate switchings, introduced in \cite{Stefanatos11} in the more general setting and confirmed in \cite{Boldt16} for the system at hand. But we also find a new kind of solution, with more than one odd number of intermediate switchings, which is absent from all the related previous works \cite{Salamon09,Tsirlin11,Salamon12,Hoffmann13,Hoffmann15,Boldt16}. This is the main contribution of this paper. Note that we didn't identify this type of solution in our previous work \cite{Stefanatos11}, in the more general case where the stiffness can take negative values, since this solution was excluded for the control set considered there.

In the next section we show the relation of the problem, as defined in \cite{Boldt16}, with our formalism in \cite{Stefanatos11}. The corresponding optimal control problem is solved in Section \ref{sec:solution}. In Section \ref{sec:examples} we illustrate the various types of solutions with several examples, and highlight the discovered new kind of solution. Section \ref{sec:conclusion} concludes the paper.

\section{Formulation of the Minimum-Time Problem}

\label{sec:formulation}

The system that we consider in this article is a particle of mass $m$ trapped in a parametric harmonic oscillator \cite{Salamon09,Tsirlin11,Salamon12,Hoffmann13,Hoffmann15,Boldt16}. The corresponding Hamiltonian is
\begin{equation}
\label{Hamiltonian}
\hat{H}=\frac{\hat{p}^2}{2m}+\frac{m\omega^2(t)\hat{q}^2}{2},
\end{equation}
where $\hat{q}, \hat{p}$ are the position and momentum operators, respectively, and $\omega(t)$ is the time-varying frequency of the oscillator which serves as the available control.
The time evolution of a quantum observable (hermitian operator) $\hat{O}$ in the Heisenberg picture is given by \cite{Merzbacher98}
\begin{equation}
\label{Observable}
\frac{d\hat{O}}{dt}=\frac{i}{\hbar}[\hat{H},\hat{O}]+\frac{\partial\hat{O}}{\partial t},
\end{equation}
where $\imath=\sqrt{-1}$ and $\hbar$ is Planck's constant.
The following operators form a closed set under the time evolution generated by $\hat{H}$ \cite{Boldt16}
\begin{equation}
\label{z_operators}
\hat{z}_1=m\hat{q}^2,\quad \hat{z}_2=\frac{\hat{p}^2}{m},\quad \hat{z}_3=-\frac{\imath}{2\hbar}[\hat{z}_1,\hat{z}_2]=\hat{q}\hat{p}+\hat{p}\hat{q}.
\end{equation}
It is sufficient to follow the expectation values
\begin{equation}
\label{expectations}
z_i=\langle\hat{z}_i\rangle=\mbox{Tr}(\rho_0\hat{z}_i),\quad i=1, 2, 3
\end{equation}
of these operators, where $\rho_0$ is the density matrix corresponding to the initial state of the system at $t=0$ (recall that we use the Heisenberg picture).
From (\ref{Observable}) and (\ref{expectations}) we easily find
\begin{align}
\label{z1}\dot{z}_1  &  = z_3,\\
\label{z2}\dot{z}_2  &  = -\omega^2z_3,\\
\label{z3}\dot{z}_3  &  = -2\omega^2z_1+2z_2.
\end{align}

In order to find the initial and final values of $z_i$ note that states of thermodynamic equilibrium, with $\omega(t)=\omega$ constant, are characterized by the equipartition of energy $E=\langle\hat{H}\rangle$
\begin{equation}
\label{equipartition}
\left\langle\frac{\hat{p}^2}{2m}\right\rangle=\left\langle\frac{m\omega^2\hat{q}^2}{2}\right\rangle=\frac{E}{2}
\end{equation}
and the absence of correlations
\begin{equation}
\label{no_correlation}
\langle\hat{q}\hat{p}+\hat{p}\hat{q}\rangle=0.
\end{equation}
If the system starts at $t=0$ from the equilibrium state with frequency $\omega_0$ and energy $E_0$, using (\ref{equipartition}) and (\ref{no_correlation}) in (\ref{z_operators}) we find
\begin{equation}
\label{z_initial}
z_1(0)=\frac{E_0}{\omega_0^2},\quad z_2(0)=E_0,\quad z_3(0)=0.
\end{equation}
For the final state at $t=T$ with frequency $\omega_f$ and energy $E_f$, the corresponding terminal conditions are
\begin{equation}
\label{z_final}
z_1(T)=\frac{E_f}{\omega_f^2},\quad z_2(T)=E_f,\quad z_3(T)=0.
\end{equation}
It can be easily verified that, during the evolution of the system, the following quantity, called the Casimir companion, is a constant of the motion \cite{Boldt13}
\begin{equation}
\label{constant}
z_1z_2-\frac{z_3^2}{4}=\frac{E_0^2}{\omega_0^2}.
\end{equation}
For all the equilibrium states it is $z_3=0$, thus these states lie on the hyperbola
\begin{equation}
\label{hyperbola}
z_1z_2=\frac{E_0^2}{\omega_0^2}
\end{equation}
in the $z_1z_2$-plane. Using (\ref{z_final}) in (\ref{hyperbola}) we find
\begin{equation}
\label{energy_ratio}
\frac{E_f}{E_0}=\frac{\omega_f}{\omega_0}.
\end{equation}
For a canonical ensemble of quantum harmonic oscillators, the equilibrium energy $E$ is related to the temperature $\mathcal{T}$ and the frequency $\omega$ through the expression
\begin{equation}
\label{temperature}
E=\frac{\hbar\omega}{2}\coth\left(\frac{\hbar\omega}{2k_b\mathcal{T}}\right).
\end{equation}
From this relation and the ratio of the energies at the initial and final states, we conclude that the corresponding ensemble temperatures satisfy
\begin{equation}
\label{temp_ratio}
\frac{\mathcal{T}_f}{\mathcal{T}_0}=\frac{\omega_f}{\omega_0}.
\end{equation}
For the case
\begin{equation}
\label{cooling}
\omega_f<\omega_0
\end{equation}
that we study here, this corresponds to a temperature reduction (cooling) by a factor of $\omega_0/\omega_f$.
We would like to find the time-varying frequency $\omega(t)$, with
\begin{equation}
\label{frequency_boundary}
\omega(t)=\left\{\begin{array}{cl} \omega_0, & t\leq 0 \\\omega_f, & t\geq T\end{array}\right.
\end{equation}
and
\begin{equation}
\label{frequency}
\omega_1\leq\omega(t)\leq\omega_2,\quad 0<t<T,
\end{equation}
where the bounds $\omega_1, \omega_2$ satisfy
\begin{equation}
\label{bounds}
0<\omega_1\leq\omega_f<\omega_0\leq\omega_2<\infty,
\end{equation}
which drives the system from the equilibrium state (\ref{z_initial}) to the equilibrium state (\ref{z_final}) in minimum time $T$.

In order to solve this problem, we will use the constant of the motion (\ref{constant}) to reduce the dimension of the system from three to two, following a different approach than that in \cite{Boldt16}. Let us define the dimensionless variable $b$ through the relations
\begin{equation}
\label{b}
b=\frac{\sqrt{\langle\hat{q}^2\rangle}}{q_0},\quad q_0=\sqrt{\frac{E_0}{m\omega_0^2}},
\end{equation}
where note that $q_0$ has length dimensions. Then, using the definition of $\hat{z}_1$ from (\ref{z_operators}) and Eqs. (\ref{z1})-(\ref{z3}), variables $z_i$ can be expressed in terms of $b$ as follows
\begin{equation}
\label{zeta}z_1=\frac{E_0}{\omega_0^2}b^2,\, z_2=\frac{E_0}{\omega_0^2}(b\ddot{b}+\dot{b}^2+\omega^2b^2),\, z_3=\frac{2E_0}{\omega_0^2}b\dot{b}.
\end{equation}
If we plug (\ref{zeta}) in (\ref{constant}), we obtain the following Ermakov equation for $b(t)$ \cite{Ermakov,Chen10}
\begin{equation}
\label{Ermakov}\ddot{b}(t)+\omega^{2}(t)b(t)=\frac{\omega_{0}^{2}}{b^{3}(t)}%
\end{equation}
The boundary conditions for $b$ can be found by using (\ref{zeta}) in (\ref{z_initial}) and (\ref{z_final}). They are
\begin{equation}
\label{b_initial}b(0)=1,\quad \dot{b}(0)=0
\end{equation}
and
\begin{equation}
\label{b_final}b(T)=\sqrt{\frac{\omega_0}{\omega_f}},\quad \dot{b}(T)=0,
\end{equation}
where we have additionally used (\ref{energy_ratio}) in the derivation of $b(T)$. Note that from (\ref{Ermakov}) we can also obtain $\ddot{b}(0)=\ddot{b}(T)=0$, as long as the frequency boundary conditions (\ref{frequency_boundary}) are satisfied, so we only need to consider the latter.

If we set
\begin{equation}
\label{x}
x_{1}=b,\quad x_{2}=\frac{\dot{b}}{\omega_{0}},\quad u(t)=\frac{\omega^{2}(t)}%
{\omega_{0}^{2}},
\end{equation}
and rescale time according to $t_{\mbox{new}}=\omega_{0} t_{\mbox{old}}$, we
obtain the following system of first order differential equations, equivalent
to the Ermakov equation 
\begin{align}
\label{system1}\dot{x}_{1}  &  = x_{2},\\
\label{system2}
\dot{x}_{2}  &  = -ux_{1}+\frac{1}{x_{1}^{3}}.
\end{align}
The control bounds are
\begin{equation}
\label{u_bounds}
u_1=\frac{\omega_1^2}{\omega_0^2}, \quad u_2=\frac{\omega_2^2}{\omega_0^2},
\end{equation}
and if we set
\begin{equation}
\label{gamma}
\gamma=\sqrt{\frac{\omega_0}{\omega_f}}>1
\end{equation}
then (\ref{bounds}) becomes
\begin{equation}
\label{u_order}
0<u_1\leq\frac{1}{\gamma^4}<1\leq u_2<\infty.
\end{equation}
Using (\ref{x}) to translate the boundary conditions (\ref{frequency_boundary}), (\ref{b_initial}), (\ref{b_final}) for $\omega, b$ into corresponding conditions for $u, x_1, x_2$, we obtain the following time-optimal problem for system (\ref{system1}), (\ref{system2}):

\newtheorem{problem}{problem} \begin{prob}\label{problem}
Find $u_1\leq u(t)\leq u_2$ with $u_1,u_2$ satisfying (\ref{u_order}) and $u(0)=1, u(T)=1/\gamma^4$, such that starting from $(x_1(0),x_2(0))=(1,0)$, the above system reaches the final point $(x_1(T),x_2(T))=(\gamma,0), \gamma>1$, in minimum time $T$.
\end{prob}



In the next section we solve the following optimal control problem, where we drop the boundary conditions on the control $u$, as we justify below:

\begin{prob}\label{problem1}
Find $u_1\leq u(t)\leq u_2$, with $u_1,u_2$ satisfying (\ref{u_order}), such that starting from $(x_1(0),x_2(0))=(1,0)$, the system above reaches the final point $(x_1(T),x_2(T))=(\gamma,0), \gamma>1$, in minimum time $T$.
\end{prob}

In both problems the class of admissible controls formally are Lebesgue measurable functions which take values in the control set $[u_1,u_2]$ almost everywhere. However, as we shall see, optimal controls are piecewise continuous, in fact bang-bang. The optimal control found for problem \ref{problem1} is also optimal for problem \ref{problem}, with the addition of instantaneous jumps at the initial and final points, so that the boundary conditions $u(0)=1$ and $u(T)=1/\gamma^4$ are satisfied. Note that in connection with (\ref{frequency_boundary}), a natural way to think about these conditions is that $u(t)=1$ for $t\leq 0$ and $u(t)=1/\gamma^4$ for $t\geq T$; in the interval $(0,T)$ we pick the control that achieves the desired transfer in minimum time.

\begin{rem}
\label{mech_analog}
Observe that the above
system (\ref{system1}), (\ref{system2}) can be interpreted as describing the one-dimensional Newtonian motion of a unit-mass particle, with position coordinate $x_{1}$ and velocity $x_{2}$. The acceleration (force) acting on the particle is $-ux_{1}+1/x_{1}^{3}$. This point of view
can provide useful intuition about the time-optimal solution, as we will see later.
\end{rem}



\section{Optimal Solution}

\label{sec:solution}

In our previous work \cite{Stefanatos11} we solved a problem similar to Problem \ref{problem1}, where the control was restricted as $-u_1\leq u(t)\leq u_2$, with $u_1,u_2\geq 1$. Note that in this setting,
the possibility of negative control $\omega^{2}(t)<0$ (expulsive parabolic potential) for some time intervals was permitted \cite{Chen10}. In the present article we consider the very interesting practical case where only attractive parabolic potential is allowed \cite{Salamon09}. In this section we investigate how our previous solution is modified due to the restriction of the control in a more narrow set. As we will see, a new type of solution arises, which was forbidden in the previous setting. In the following, we provide the details of the proofs which are modified, compared to the previous case, and also the basic steps of the proofs which remain the same, for completeness.

The system described by (\ref{system1}), (\ref{system2}) can be expressed
in compact form as
\begin{equation}
\dot{x}=f(x)+ug(x), \label{affine}%
\end{equation}
where the vector fields are given by
\begin{equation}
f=\left(
\begin{array}
[c]{c}%
x_{2}\\
1/x_{1}^{3}%
\end{array}
\right)  ,\,\,g=\left(
\begin{array}
[c]{c}%
0\\
-x_{1}%
\end{array}
\right)
\end{equation}
and $x\in\mathcal{D}=\{(x_{1},x_{2})\in\mathbb{R}^{2}:x_{1}>0\}$ and $u\in
U=[u_{1},u_{2}]$. Admissible controls are Lebesgue measurable functions that
take values in the control set $U$. Given an admissible control $u$ defined
over an interval $[0,T]$, the solution $x$ of the system (\ref{affine})
corresponding to the control $u$ is called the corresponding trajectory and we
call the pair $(x,u)$ a controlled trajectory. Note that the domain
$\mathcal{D}$\ is invariant in the sense that trajectories cannot leave
$\mathcal{D}$. Starting with any positive initial condition $x_{1}(0)>0$, and
using any admissible control $u$, as $x_{1}\rightarrow0^{+}$, the
``repulsive force" $1/x_{1}^{3}$ leads to an increase in
$x_{1}$ that will keep $x_{1}$ positive (as long as the solutions exist).

For a constant $\lambda_{0}$ and a row vector $\lambda=(\lambda_{1}%
,\lambda_{2})\in\left(  \mathbb{R}^{2}\right)  ^{\ast}$ define the control
Hamiltonian as%
\[
H=H(\lambda_{0},\lambda,x,u)=\lambda_{0}+\langle\lambda,f(x)+ug(x)\rangle.
\]
Pontryagin's Maximum Principle for \emph{time-optimal} processes \cite{Pontryagin}
provides the following necessary conditions for optimality:

\begin{theorem}
[\textrm{Maximum principle}%
]\cite{Pontryagin} \label{prop:max_principle} Let $(x_{\ast}(t),u_{\ast}(t))$
be a time-optimal controlled trajectory that transfers the initial condition
$x(0)=x_{0}$ into the terminal state $x(T)=x_T$. Then it is a necessary
condition for optimality that there exists a constant $\lambda_{0}\leq0$ and
nonzero, absolutely continuous row vector function $\lambda(t)$ such that:

\begin{enumerate}
\item $\lambda$ satisfies the so-called adjoint equation%
\[
\dot{\lambda}(t)=-\frac{\partial H}{\partial x}(\lambda_{0},\lambda
(t),x_{\ast}(t),u_{\ast}(t))
\]

\item For $0\leq t\leq T$ the function $u\mapsto H(\lambda_{0}%
,\lambda(t),x_{\ast}(t),u)$ attains its maximum\ over the control set $U$ at
$u=u_{\ast}(t)$.

\item $H(\lambda_{0},\lambda(t),x_{\ast}(t),u_{\ast}(t))\equiv0$.
\end{enumerate}
\end{theorem}

We call a controlled trajectory $(x,u)$ for which there exist multipliers
$\lambda_{0}$ and $\lambda(t)$ such that these conditions are satisfied an
extremal. Extremals for which $\lambda_{0}=0$ are called abnormal. If
$\lambda_{0}<0$, then without loss of generality we may rescale the $\lambda
$'s and set $\lambda_{0}=-1$. Such an extremal is called normal.

For the system (\ref{system1}), (\ref{system2}) we have
\begin{equation}
H(\lambda_{0},\lambda,x,u)=\lambda_{0}+\lambda_{1}x_{2}+\lambda_{2}\left(  \frac{1}%
{x_{1}^{3}}-x_{1}u\right)  ,\label{hamiltonian}%
\end{equation}
and thus
\begin{equation}
\dot{\lambda}=-\lambda\left(
\begin{array}
[c]{cc}%
0 & 1\\
-(u+3/x_{1}^{4}) & 0
\end{array}
\right)  =-\lambda A\label{adjoint}%
\end{equation}

Observe that $H$ is a linear function of the bounded control variable $u$. The
coefficient at $u$ in $H$ is $-\lambda_{2}x_{1}$ and, since $x_{1}>0$, its
sign is determined by $\Phi=-\lambda_{2}$, the so-called \emph{switching
function}. According to the maximum principle, point 2 above, the optimal
control is given by $u=u_{1}$ if $\Phi<0$ and by $u=u_{2}$ if $\Phi>0$. The
maximum principle provides a priori no information about the control at times
$t$ when the switching function $\Phi$ vanishes. However, if $\Phi(t)=0$ and
$\dot{\Phi}(t)\neq0$, then at time $t$ the control switches between its
boundary values and we call this a bang-bang switch. If $\Phi$ were to vanish
identically over some open time interval $I$ the corresponding control is
called \emph{singular}.

\begin{prop}
For Problem \ref{problem1} optimal controls are bang-bang.
\end{prop}

\begin{IEEEproof}
Whenever the switching function $\Phi(t)=-\lambda_{2}(t)$ vanishes at some
time $t$, then it follows from the non-triviality of the multiplier
$\lambda(t)$ that its derivative $\dot{\Phi}(t)=-\dot{\lambda}_{2}%
(t)=\lambda_{1}(t)$ is non-zero. Hence the switching function changes sign and
there is a bang-bang switch at time $t$.
\end{IEEEproof}


\begin{mydef}
We denote the vector fields corresponding to the constant bang controls
$u_{1}$ and $u_{2}$ by $X=f+u_{1}g$ and $Y=f+u_{2}g$, respectively, and call
the trajectories corresponding to the constant controls $u\equiv u_{1}$ and
$u\equiv u_{2}$ $X$- and $Y$-trajectories. A concatenation of an
$X$-trajectory followed by a $Y$-trajectory is denoted by $XY$ while the
concatenation in the inverse order is denoted by $YX$.
\end{mydef}

We next show that all the extremals of the problem are normal. We use the following lemma:

\begin{lem}
\label{abnormal}
An $X$-trajectory starting from $(\alpha,0), 0<\alpha\leq 1$, meets the $x_1$-axis at a point $(\beta,0)$ with $\beta>\gamma$. A $Y$-trajectory starting from $(\beta,0), \beta\geq 1$, meets the $x_1$-axis at a point $(\alpha,0)$ with $\alpha\leq1$.
\end{lem}

\begin{IEEEproof}
For an $X$-trajectory ($u=u_1$) starting from $(\alpha,0)$ it is not hard to verify, using the system equations, the following constant of the motion
\begin{equation}
x_{2}^{2}+u_{1}x_{1}^{2}+\frac{1}{x_{1}^{2}}=u_{1}\alpha^{2}+\frac{1}{\alpha^{2}}.\label{first_integral_1}
\end{equation}
For $x_2=0$ the above equation has two solutions for $x_1$: $\alpha$, corresponding to the starting point, and $\beta=1/(\alpha\sqrt{u_1})$. But $u_1\leq 1/\gamma^4$ from (\ref{u_order}) and $\alpha\leq 1$ from lemma hypothesis, thus $\beta\geq\gamma^2>\gamma$, since additionally $\gamma>1$. Analogously, a first integral of the motion along the $Y$-trajectory ($u=u_2$) starting from $(\beta,0)$ is
\begin{equation}
x_{2}^{2}+u_{2}x_{1}^{2}+\frac{1}{x_{1}^{2}}=u_{2}\beta^{2}+\frac{1}{\beta^{2}}.\label{first_integral_2}
\end{equation}
For $x_2=0$ we obtain two values for $x_1$, $\beta$ (starting point) and $\alpha=1/(\beta\sqrt{u_2})$. But $\beta\geq1$ and $u_2\geq 1$, thus $\alpha\leq1$.
\end{IEEEproof}

\begin{prop}
All the extremals are normal.
\end{prop}

\begin{IEEEproof}
If $(x,u)$ is an abnormal extremal trajectory with a switching at
$t=t_0$, then, since $\lambda_{2}(t_0)=0$, it follows from $H=0$ that it is also
$x_{2}(t_0)=0$. Thus, for abnormal extremals, all the switchings take place on the $x_1$-axis.
Suppose now that the system starts from $(\alpha_1=1,0)$ with an $X$-segment ($u=u_1$). This trajectory meets again the $x_1$-axis at a point $(\beta_1,0)$, where $\beta_1>\gamma>1$, according to the above lemma. At this point there is a switching to $u=u_2$, otherwise the system returns to the starting point. The $Y$-segment starting from $(\beta_1,0)$, $\beta_1>1$, meets again the $x_1$-axis at a point $(\alpha_2,0)$, where $\alpha_2<1$, according to the lemma. By repeating this procedure, we observe that the abnormal extremal trajectory is passing from $x_1$-axis only through points $(\alpha_i,0), (\beta_i,0)$, with $\alpha_i<1$ and $\beta_i>\gamma$, thus it can never reach the target point $(\gamma,0)$. The proof is analogous when the trajectory starts with a $Y$-segment from $(\beta_1=1,0)$.
\end{IEEEproof}

For normal extremals we can set $\lambda_{0}=-1$.
Then, $H=0$ implies that for any switching time $t_0$ we must have $\lambda
_{1}(t_0)x_{2}(t_0)=1$. For an $XY$ junction we have $\dot{\Phi}%
(t_0)=\lambda_{1}(t_0)>0$ and thus necessarily $x_{2}(t_0)>0$ and analogously
optimal $YX$ junctions need to lie in $\{x_{2}<0\}$. In the following, we establish the precise concatenation sequences for optimal
controls and in particular calculate the times between switchings explicitly.

\begin{lem}
[\textrm{Inter-switching time}]\label{switch_time} Let $P=(x_{1},x_{2})$ be a
switching point and $\tau$ denote the time to reach the next switching point
$Q$. If $\overrightarrow{PQ}$ is a $Y$-trajectory, then
\begin{equation}
\sin(2\sqrt{u_{2}}\tau)=-\frac{2\sqrt{u_{2}}x_{1}x_{2}}{x_{2}^{2}+u_{2}%
x_{1}^{2}},\quad\cos(2\sqrt{u_{2}}\tau)=\frac{x_{2}^{2}-u_{2}x_{1}^{2}}%
{x_{2}^{2}+u_{2}x_{1}^{2}}\label{par_2}%
\end{equation}
while, if $\overrightarrow{PQ}$ is an $X$-trajectory, then
\begin{equation}
\sin(2\sqrt{u_{1}}\tau)=-\frac{2\sqrt{u_{1}}x_{1}x_{2}}{x_{2}^{2}+u_{1}%
x_{1}^{2}},\quad\cos(2\sqrt{u_{1}}\tau)=\frac{x_{2}^{2}-u_{1}x_{1}^{2}}%
{x_{2}^{2}+u_{1}x_{1}^{2}}.\label{par_1}%
\end{equation}

\end{lem}

Note that the inter-switching times depend only on the ratio $x_2/x_1$.

\begin{IEEEproof}
These formulas are obtained as an application of the concept of a
``conjugate point" for bang-bang controls \cite{Suplane,Boscain04}. The proof is the
same as in \cite{Stefanatos11}, but for completeness we repeat here the main steps.
Without loss of generality assume that the trajectory passes through $P$ at
time $0$ and is at $Q$ at time $\tau$. Since $P$ and $Q$ are switching points,
the corresponding multipliers vanish against the control vector field $g$ at
those points, i.e., $\langle\lambda(0),g(P)\rangle=\langle\lambda
(\tau),g(Q)\rangle=0$. We need to compute what the relation $\langle
\lambda(\tau),g(Q)\rangle=0$ implies at time $0$. In order to do so, we move
the vector $g(Q)$ along the $Y$-trajectory backward from $Q$ to $P$. This is
done by means of the solution $w(t)$ of the variational equation along the
$Y$-trajectory with terminal condition $w(\tau)=g(Q)$ at time $\tau$. Recall
that the variational equation along $Y$ is the linear system $\dot{w}=Aw$
where matrix $A$ is given in (\ref{adjoint}). Symbolically, if we denote by
$e^{tY}(P)$ the value of the $Y$-trajectory at time $t$ that starts at the
point $P$ at time $0$ and by $(e^{-tY})_{\ast}$ the backward evolution under
the linear differential equation $\dot{w}=Aw$, then we can represent this
solution in the form
\begin{align}
w(0)&=(e^{-\tau Y})_{\ast}w(\tau)=(e^{-\tau Y})_{\ast}g(Q)\nonumber\\&=(e^{-\tau Y})_{\ast}g(e^{\tau
Y}(P))=(e^{-\tau Y})_{\ast}\circ
g\circ e^{\tau Y}(P).\nonumber
\end{align}
Since the
\textquotedblleft adjoint equation\textquotedblright\ of the Maximum Principle
is precisely the adjoint equation to the variational equation, it follows that
the function $t\mapsto\langle\lambda(t),w(t)\rangle$ is constant along the
$Y$-trajectory. Hence $\langle\lambda(\tau),g(Q)\rangle=0$ implies that
\[
\langle\lambda(0),w(0)\rangle=\langle\lambda(0),(e^{-\tau Y})_{\ast}g(e^{\tau
Y}(P))\rangle=0
\]
as well. But the non-zero two-dimensional multiplier $\lambda(0)$ can only be orthogonal to
both $g(P)$ and $w(0)$ if these vectors are parallel, $g(P)\Vert
w(0)=(e^{-\tau Y})_{\ast}g(e^{\tau Y}(P))$. It is this relation that defines
the switching time.

It remains to compute $w(0)$. For this we make use of the well-known relation \cite{Heinz12}
\begin{equation}
(e^{-\tau Y})_{\ast}\circ g\circ e^{\tau Y}=e^{\tau\,adY}(g)\label{adj}%
\end{equation}
where the operator $adY$ is defined as $adY(g)=[Y,g]$, with $[,]$ denoting the
Lie bracket of the vector fields $Y$ and $g$. 
For our system, the Lie algebra 
generated by the fields $f$ and $g$ actually is finite dimensional: we have
\[
\lbrack f,g](x)=\left(
\begin{array}
[c]{c}%
x_{1}\\
-x_{2}%
\end{array}
\right)
\]
and the relations
\[
\lbrack f,[f,g]]=2f,\qquad\lbrack g,[f,g]]=-2g
\]
can be directly verified. Using these relations and the analyticity of the
system, $e^{t\,adY}(g)$ can be calculated in closed form from the expansion
\begin{equation}
e^{t\,adY}(g)=\sum_{n=0}^{\infty}\frac{t^{n}}{n!}\,ad\,^{n}Y(g),
\end{equation}
where, inductively, $ad^{n}Y(g)=[Y,ad^{n-1}Y(g)]$.
By summing the series appropriately we obtain
\begin{align}
e^{t\,adY}(g)&=g+\frac{1}{2\sqrt{u_{2}}}\sin(2\sqrt{u_{2}}t)[f,g]\nonumber\\&+\frac
{1}{2u_{2}}[1-\cos(2\sqrt{u_{2}}t)](f-u_{2}g).\nonumber
\end{align}
The field $w(0)=(e^{-\tau Y})_{\ast}g(e^{\tau Y}(P))=e^{t\,adY}(g(P))$ is
parallel to $g(P)=(0,-x_{1})^{T}$ if and only if
\begin{equation}
\sqrt{u_{2}}x_{1}\sin(2\sqrt{u_{2}}\tau)+x_{2}\left[  1-\cos(2\sqrt{u_{2}}%
\tau)\right]  =0.\nonumber
\end{equation}
Hence
\begin{equation}
\sin(2\sqrt{u_{2}}\tau)=-\frac{x_{2}}{\sqrt{u_{2}}x_{1}}[1-\cos(2\sqrt{u_{2}%
}\tau)]
\end{equation}
from which (\ref{par_2}) follows. 

For an $X$-trajectory  we simply replace $u_2$ with $u_1$ and obtain (\ref{par_1}).
\end{IEEEproof}

\begin{figure}[t]
\centering
\includegraphics[width=0.8\linewidth]{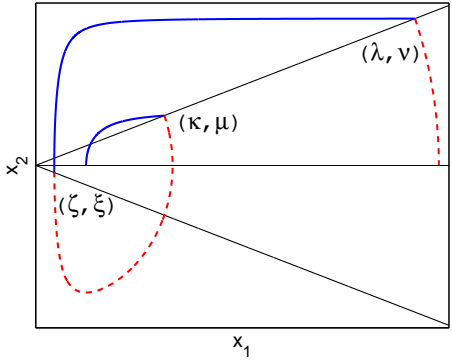}\caption{Consecutive switching points lie on two opposite-slope lines through the origin. Blue solid curves correspond to $X$-segments $(u=u_1)$, red dashed curves to $Y$-segments $(u=u_2)$.}%
\label{fig:threeswitchings}%
\end{figure}

\begin{lem}
[\textrm{Main technical point}]\label{prop:ratio} The ratio of the coordinates
of consecutive switching points has constant magnitude but alternating sign,
while these points are not symmetric with respect to the $x_{1}$-axis.
\end{lem}

\begin{IEEEproof}
Consider the trajectory $XYXY$ shown in Fig.\ \ref{fig:threeswitchings} with
switching points $(\kappa,\mu), (\zeta,\xi)$ and $(\lambda,\nu)$, where blue solid curves correspond to $X$-segments $(u=u_1)$ and red dashed curves to $Y$-segments $(u=u_2)$.
If we follow \cite{Stefanatos11} then, starting from $(\kappa,\mu)$ and integrating the equations of motion (\ref{system1}) and (\ref{system2}) for the inter-switching time given in (\ref{par_2}), we can find the coordinates of the next switching point and show that $\xi/\zeta=-\mu/\kappa$ while $(\zeta,\xi)\neq (\kappa,-\mu)$. Subsequently, integrating the equations for the inter-switching time given in (\ref{par_1}), we can also show that $\nu/\lambda=-\xi/\zeta$ and $(\lambda,\nu)\neq (\zeta,-\xi)$.

Here we present a more elegant proof based on the symmetries of the system. Observe that the transformation $(t,x_1,x_2)\rightarrow (-t,x_1,-x_2)$ leaves the system (\ref{system1}) and (\ref{system2}) invariant for constant $u$. So, starting from $(\zeta,-\xi)$ and running the \emph{transformed} system forward in time, we arrive at the next switching point $(\kappa,\mu)$, which is the point encountered when running the original system backward in time. The switching time is given again by (\ref{par_2}), with $x_1=\zeta, x_2=-\xi$. But this switching time is the same when running the original system in the forward direction, with $x_1=\kappa, x_2=\mu$ in (\ref{par_2}). Equating the sine and cosine terms in (\ref{par_2}) for the forward and backward directions, we obtain
\begin{equation}
(u_2\zeta\kappa+\xi\mu)(\xi\kappa+\zeta\mu)=0,\quad(\xi\kappa-\zeta\mu)(\xi\kappa+\zeta\mu)=0.
\end{equation}
The above equations are both satisfied when $\xi/\zeta=-\mu/\kappa$. Note that it is $(\zeta,\xi)\neq (\kappa,-\mu)$, otherwise the trajectory would return to the starting point shown in Fig. \ref{fig:threeswitchings} forming a loop, which is obviously not time-optimal.
\end{IEEEproof}

\begin{lem}
\label{integral}
Starting from the point $(\beta,0)$ at $t=0$, the time evolution of $x_1$ is
\begin{equation}
\label{x1_in_Y}x_{1}^2(t)=\frac{1}{2}\left(  \beta^{2}+\frac{1}{u%
\beta^{2}}\right)  +\frac{1}{2}\left(  \beta^{2}-\frac{1}{u\beta^{2}%
}\right)  \cos(2\sqrt{u}t),
\end{equation}
where $u=u_1$ for an $X$-segment and $u=u_2$ for a $Y$-segment.
\end{lem}

\begin{IEEEproof}
Integrate the system equations (\ref{system1}), (\ref{system2}) as in \cite{Stefanatos11}.
\end{IEEEproof}
In the following Theorem \ref{prop:time}, which is the main result of the paper, we use Lemma \ref{prop:ratio} to determine the exact form
of the extremals, and Lemmas \ref{switch_time} and \ref{integral} to calculate the corresponding times to reach the target point.




\begin{theorem}
\label{prop:time}
The extremal trajectories can only have the form $XYX\ldots XY$, with an odd number of switchings, or the form
$YXY\ldots XY$, with an even number of switchings. The necessary time to reach the
target point $(\gamma,0), \gamma>1$, with the extremal $XYX\ldots XY$ with $2n+1$ switchings, $n=0,1,2,\ldots$, is
\begin{equation}
\label{time_odd}
T^{\pm}_{2n+1}=T^{\pm}_{I,1}+n(T_{X}+T_{Y})+T_{F},
\end{equation}
where
\begin{align}
\label{time_in1}T^{\pm}_{I,1}  &  = \frac{1}{2\sqrt{u_{1}}}\cos^{-1}\left(
\frac{sc_{1}\mp u_{1}\sqrt{c_{1}^{2}-4(s+u_{1})}}{(s+u_{1})\sqrt{c_{1}%
^{2}-4u_{1}}}\right),\\
\label{time_fi}
T_{F}  &  = \frac{1}{2\sqrt{u_{2}}}\cos^{-1}\left(  \frac{-sc+u_{2}%
\sqrt{c^{2}-4(s+u_{2})}}{(s+u_{2})\sqrt{c^{2}-4u_{2}}}\right)  ,
\end{align}
\begin{align}
\label{switch_X}
T_{X}  &  = \frac{1}{2\sqrt{u_{1}}}\cos^{-1}\left(
\frac{s-u_{1}}{s+u_{1}}\right),\\
\label{switch_Y}
T_{Y}  &  = \frac{1}{2\sqrt{u_{2}}}\left(  2\pi-\cos^{-1}\left(  \frac
{s-u_{2}}{s+u_{2}}\right)  \right)  ,
\end{align}
\begin{align}
\label{c_1}
c_{1}  &  =u_{1}+1,\\
\label{c}
c  &  =u_{2}\gamma^{2}+\frac{1}{\gamma^{2}},
\end{align}
and $s$ is the solution of the transcendental equation
\begin{equation}
\label{transcendentalX}\frac{c+\sqrt{c^{2}-4(s+u_{2})}}{c_{1}%
\pm\sqrt{c_{1}^{2}-4(s+u_{1})}}=\left(\frac{s+u_{2}}{s+u_{1}}\right)^{n+1}%
\end{equation}
in the interval $0<s\leq\mbox{Min}\{(1-u_{1})^{2}/4, (u_{2}\gamma^2-1/\gamma^2)^{2}/4\}$. Note that the $\pm$ sign in (\ref{transcendentalX}) corresponds to the $\pm$ sign in (\ref{time_odd}). The constants $c_{1}$ and
$c$ characterize the first $X$-segment and the last $Y$-segment, respectively, of
the trajectory. The necessary time to reach the
target point with the extremal $YXY\ldots XY$ with $2n$ switchings, $n=1,2,\ldots$, is
\begin{equation}
\label{time_even}
T^{\pm}_{2n}=T^{\pm}_{I,2}+nT_{X}+(n-1)T_{Y}+T_{F},
\end{equation}
where
\begin{equation}
\label{time_in2}T^{\pm}_{I,2} = \frac{1}{2\sqrt{u_{2}}}\cos^{-1}\left(
\frac{-sc_{2}\pm u_{2}\sqrt{c_{2}^{2}-4(s+u_{2})}}{(s+u_{2})\sqrt{c_{2}%
^{2}-4u_{2}}}\right),
\end{equation}
$T_X, T_Y, T_F$ are the same as above,
\begin{equation}
\label{c_2}
c_{2}=u_{2}+1,
\end{equation}
and $s$ is the solution of the transcendental equation
\begin{equation}
\label{transcendentalY}\frac{c+\sqrt{c^{2}-4(s+u_{2})}}{c_{2}%
\mp\sqrt{c_{2}^{2}-4(s+u_{2})}}=\left(\frac{s+u_{2}}{s+u_{1}}\right)^{n}%
\end{equation}
in the interval $0<s\leq(u_{2}-1)^{2}/4$. The $\mp$ sign in (\ref{transcendentalY}) corresponds to the $\pm$ sign in (\ref{time_even}), while the constant (\ref{c_2}) characterizes the first $Y$-segment of the trajectory.
\end{theorem}

\begin{IEEEproof}
Consider a trajectory of the form $XYX\ldots XY$ with $n$ turns and $2n+1$ switching points $A_j(\kappa_{j}%
,\mu_{j}), j=1,2,\ldots,2n+1$, shown in Fig. \ref{fig:odd}. Observe that the odd-numbered switching points lie on a positive-slope straight line passing through the origin, while the even-numbered switching points lie on the symmetric line with opposite slope, in accordance with Lemma \ref{prop:ratio}. Two consecutive switching points satisfy the following equation
\begin{equation}
\label{connection}
\mu_{j+1}^{2}+u\kappa_{j+1}^{2}+\frac{1}{\kappa_{j+1}^{2}}=\mu_{j}^{2}+u\kappa_{j}^{2}+\frac{1}{\kappa_{j}^{2}},
\end{equation}
where $u=u_1$ if the two points are connected with an $X$-segment and $u=u_2$ if they are joined with a $Y$-segment (it can be verified from the system equations that the quantity $x_2^2+ux_1^2+1/x_1^2$ is constant along segments with constant control $u$).
The \emph{ratio of the squares} of the coordinates of all the switching points is constant, and if we denote it by $\mu_{j+1}^{2}/\kappa_{j+1}^{2}=\mu_{j}^{2}/\kappa_{j}^{2}=s$, then (\ref{connection}) becomes
\begin{equation}
(\kappa_{j+1}^{2}-\kappa_{j}^{2})\left(s+u-\frac{1}{\kappa_{j}^{2}\kappa_{j+1}^{2}}\right)=0
\end{equation}
But $\kappa_{j+1}\neq\kappa_{j}$ since the consecutive switching points are not symmetric with respect to $x_1$-axis (Lemma \ref{prop:ratio}), thus
\begin{equation}
\label{consecutive}
\kappa_{j+1}^{2}=\frac{1}{\kappa_{j}^{2}(s+u)}.
\end{equation}
If we apply (\ref{consecutive}) for three successive switching points we obtain
\begin{equation}
\label{X_odd}
\frac{\kappa_{2k+1}^{2}}{\kappa_{2k-1}^{2}}=\frac{s+u_2}{s+u_1}>1\Rightarrow \kappa_{1}<\kappa_{3}<\ldots<\kappa_{2n+1}
\end{equation}
for the odd switching points and
\begin{equation}
\label{X_even}
\frac{\kappa_{2k+2}^{2}}{\kappa_{2k}^{2}}=\frac{s+u_1}{s+u_2}<1\Rightarrow \kappa_{2}>\kappa_{4}>\ldots>\kappa_{2n}
\end{equation}
for the even switching points. We show that an extremal starting with an $X$-segment cannot also end with an $X$-segment. If that was true, then the last switching point would have even numbering $(\kappa_{2k},\mu_{2k})$, leading to a final point on the $x_1$-axis with $x_1(T)<\kappa_{2k}$. This happens because at this switching point it is $x_2<0$, a negative velocity according to the particle model from Remark \ref{mech_analog}, thus the state of the system moves to smaller $x_1$ for the repulsive force $1/x_1^3$ in (\ref{system2}) to reduce the magnitude of the final velocity to zero $x_2(T)=0$. From the ordering in (\ref{X_even}) we conclude that it is also $x_1(T)<\kappa_{2}$. But the first two switching points satisfy (\ref{connection}) with $u=u_2$, since they are connected with a $Y$-segment, and if we use the common ratio $\mu_{2}^{2}/\kappa_{2}^{2}=\mu_{1}^{2}/\kappa_{1}^{2}=s$, we easily find that both $\kappa_1^2,\kappa_2^2$ are the roots of the following equation
\begin{equation}
\label{trionymo}
(s+u_2)\kappa^4-C\kappa^2+1=0,
\end{equation}
where $C=\mu_1^2+u_2\kappa_1^2+1/\kappa_1^2=\mu_2^2+u_2\kappa_2^2+1/\kappa_2^2$. Thus
\begin{equation}
\label{kappa12}
\kappa_{1}^2\kappa_{2}^2=\frac{1}{s+u_2}<1,
\end{equation}
since $u_2\geq 1$ and $s>0$. The first switching point belongs to an $X$-segment starting from $(1,0)$, and one can easily show that
\begin{equation}
\label{kappa1}
1\leq\kappa_1\leq 1/\sqrt{u_1},
\end{equation}
where recall that $u_1<1$. From (\ref{kappa12}), (\ref{kappa1}) we conclude that $\kappa_2<1$, thus also $x_1(T)<\kappa_2<1<\gamma$, and the final point $(\gamma,0)$ cannot be reached. Consequently, an extremal starting with an $X$-segment can only end with a $Y$-segment, as we considered at the beginning.

\begin{figure}[t]
\centering
\begin{tabular}
[c]{c}%
\subfigure[$\ $XYX\ldots XY$$]{ \label{fig:odd} \includegraphics[width=0.9\linewidth]{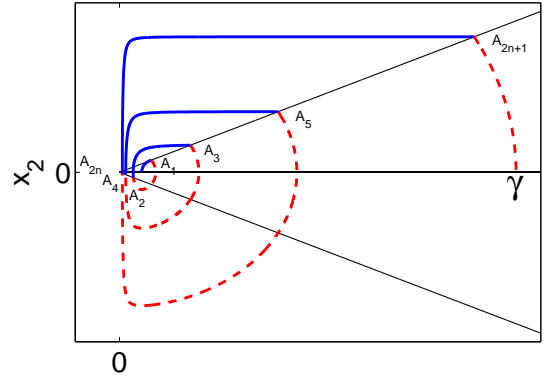}}\\
\subfigure[$\ $YXY\ldots XY$$]{ \label{fig:even} \includegraphics[width=0.9\linewidth]{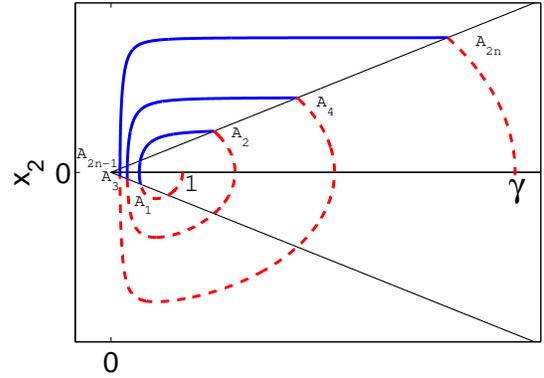}}
\end{tabular}
\caption{The two types of extremals, where blue solid line corresponds to $X$-segments ($u=u_1$) and red dashed line corresponds to $Y$-segments ($u=u_2$): (a) $XYX\ldots XY$ with odd number of switchings (b) $YXY\ldots XY$ with even number of switchings.}%
\label{fig:proof}%
\end{figure}

We next move to find an equation for the ratio $s$. If we consecutively apply (\ref{consecutive}) from the first switching point up to the last, we obtain
\begin{equation}
\label{X_first_last}
\frac{\kappa_{2n+1}^2}{\kappa_{1}^2}=\left(\frac{s+u_2}{s+u_1}\right)^n.
\end{equation}
Since the first switching point belongs to the first $X$-segment starting from $(1,0)$, it satisfies the equation
\begin{equation}
\label{k1_trionymo}
(s+u_1)\kappa_1^4-c_1\kappa_1^2+1=0,
\end{equation}
where $c_1=u_1+1$. Solving for $\kappa_1^2$ we obtain
\begin{equation}
\label{k1}
\kappa_{1,\pm}^{2}=\frac{c_{1}\pm\sqrt{c_{1}^{2}-4(s+u_{1})}}{2(s+u_{1})}=\frac{2}{c_{1}\mp\sqrt{c_{1}^{2}-4(s+u_{1})}}.
\end{equation}
The last switching point $A_{2n+1}$ belongs to the $X$-segment just before the last $Y$-segment, and satisfies an equation of the form
\begin{equation}
\label{k_last_trionymo}
(s+u_1)\kappa_{2n+1}^4-\bar{c}\kappa_{2n+1}^2+1=0,
\end{equation}
where $\bar{c}=u_1\beta^2+1/\beta^2$ and $(\beta,0)$ is the point where this $X$-segment would meet the $x_1$-axis if continued.
Note that this $X$-segment already meets the $x_1$-axis at another point $(\alpha,0)$ with $\alpha\leq1$ (in Fig. \ref{fig:odd} it is the trajectory point closest to the origin), as explained in the previous paragraph, thus, from the proof of Lemma \ref{abnormal} it is $\beta= 1/(\alpha\sqrt{u_1})\geq 1/\sqrt{u_1}>1$.
Solving (\ref{k_last_trionymo}) we find
\begin{equation}
\label{k_last}
\kappa_{2n+1,\pm}^{2}=\frac{\bar{c}\pm\sqrt{\bar{c}^{2}-4(s+u_{1})}}{2(s+u_{1})}=\frac{2}{\bar{c}\mp\sqrt{\bar{c}^{2}-4(s+u_{1})}}.
\end{equation}
If we choose $\kappa_{2n+1,+}^{2}$, which corresponds to the $+$ sign in the first equality in (\ref{k_last}), then (\ref{X_first_last}) becomes
\begin{equation}
\label{transcendentalXplus}\frac{\bar{c}+\sqrt{\bar{c}^{2}-4(s+u_{1})}}{c_{1}%
\pm\sqrt{c_{1}^{2}-4(s+u_{1})}}=\left(\frac{s+u_{2}}{s+u_{1}}\right)^{n},%
\end{equation}
with the $\pm$ sign corresponding to $\kappa_{1,\pm}^{2}$. If we choose $\kappa_{2n+1,-}^{2}$, corresponding to the $+$ sign in the second equality in (\ref{k_last}), and use for $\kappa_{1,\pm}^{2}$ the second equality in (\ref{k1}), we obtain an equation similar to (\ref{transcendentalXplus}) but with inverted left hand side. It is $\bar{c}%
>c_{1}\Leftrightarrow(\beta^{2}-1)(u_{1}\beta^{2}-1)>0$, which is true, and $c_1,\bar{c}>0$, so
\begin{equation}
\frac{c_{1}\mp\sqrt{c_{1}^{2}-4(s+u_{1})}}{\bar{c}+\sqrt{\bar{c}^{2}-4(s+u_{1}%
)}}<1<\left(  \frac{s+u_{2}}{s+u_{1}}\right)  ^{n},\nonumber
\end{equation}
and the corresponding transcendental equation has no solution. Since we actually do not know $\bar{c}$ (we do not know $\beta$) in the valid equation (\ref{transcendentalXplus}), we will use the fact that the last switching point $(\kappa_{2n+1},\mu_{2n+1})$ belongs to the final $Y$-segment passing from the target point $(\gamma,0)$. Then, $\kappa_{2n+1}^{2}>1$ is the larger root of the equation
\begin{equation}
\label{k_last_trionymo_Y}
(s+u_2)\kappa_{2n+1}^4-c\kappa_{2n+1}^2+1=0,
\end{equation}
where $c=u_2\gamma^2+1/\gamma^2$, thus
\begin{equation}
\label{k_last_Y}
\kappa_{2n+1}^{2}=\frac{c+\sqrt{c^{2}-4(s+u_{2})}}{2(s+u_{2})}.
\end{equation}
Using (\ref{k_last_Y}) and the first equality of (\ref{k1}) in (\ref{X_first_last}), we obtain the transcendental equation (\ref{transcendentalX}) for the ratio $s$, where $c_1, c$ are given in (\ref{c_1}) and (\ref{c}), respectively. Note that $s$ is bounded below by the requirement $s>0$ (for $s=0$ the switching points would lie on the $x_1$-axis, which is not the case), and above by the requirements $c_{1}^{2}-4(s+u_{1})\geq0$ and $c^{2}-4(s+u_{2})\geq0$, which are both satisfied for $s\leq\mbox{Min}\{(1-u_1)^{2}/4, (u_{2}\gamma^2-1/\gamma^2)^{2}/4\}$. 

Once we have found this ratio, we can calculate the time interval between consecutive
switchings using (\ref{switch_X}) for an $X$-segment and (\ref{switch_Y}) for
a $Y$-segment, relations obtained from Lemma \ref{switch_time} on the inter-switching time. The difference in the two expressions comes from the fact that the sines in (\ref{par_2}), (\ref{par_1}) have opposite signs and in (\ref{switch_X}), (\ref{switch_Y}) we use the inverse cosine function with range $[0,\pi]$. Observe that the times along all intermediate $X$- (respectively $Y$-) trajectories are equal. The initial time interval $T_{I,1}^{\pm}$ from the starting point $(1,0)$ up to the first switching $A_1$ can be calculated by setting $\beta=1, u=u_1$ and $x_1(T_{I,1}^{\pm})=\kappa_{1,\pm}$ in (\ref{x1_in_Y}). The result is given in (\ref{time_in1}). Analogously, the final time interval $T_{F}$ from the last
switching $A_{2n+1}$ up to the target point $(\gamma, 0)$ can be calculated by setting $\beta=\gamma, u=u_2$ and $x_1(T_{F})=\kappa_{2n+1}$ in (\ref{x1_in_Y}), and the result is given in (\ref{time_fi}). The total duration $T_{2n+1}^{\pm}$ of
the trajectory with $2n+1$ switchings joining the points $(1,0)$ and $(\gamma,0)$ is
given by (\ref{time_odd}), where $\pm$ corresponds to the choice of sign in (\ref{k1}) for the first switching point.

Consider now an extremal of the form $YX\ldots YXY$ with $n$ turns and $2n$ switching points $(\kappa_{j}%
,\mu_{j}), j=1,2,\ldots,2n$, shown in Fig. \ref{fig:even}.
If we follow a procedure similar to the one above we find
\begin{equation}
\label{Y_odd}
\frac{\kappa_{2k+1}^{2}}{\kappa_{2k-1}^{2}}=\frac{s+u_1}{s+u_2}<1\Rightarrow \kappa_{1}>\kappa_{3}>\ldots>\kappa_{2n-1}
\end{equation}
for the odd switching points and
\begin{equation}
\label{Y_even}
\frac{\kappa_{2k+2}^{2}}{\kappa_{2k}^{2}}=\frac{s+u_2}{s+u_1}>1\Rightarrow \kappa_{2}<\kappa_{4}<\ldots<\kappa_{2n}
\end{equation}
for the even switching points. We next show that a trajectory starting with a $Y$-segment cannot end with an $X$-segment. If that was the case, the last switching point would have odd numbering $(\kappa_{2k+1},\mu_{2k+1})$, leading to a final point on the $x_1$-axis with $x_1(T)<\kappa_{2k+1}$, since at the switching point $x_2<0$ and the state of the system moves to smaller $x_1$ for the repulsive force $1/x_1^3$ to reduce the magnitude of the velocity to zero $x_2(T)=0$. From the ordering in (\ref{Y_odd}) we conclude that it is also $x_1(T)<\kappa_{1}$. The first switching point belongs to a $Y$-segment starting from $(1,0)$, and one can easily show that
\begin{equation}
\label{Y_kappa1}
1/\sqrt{u_2}\leq\kappa_1\leq 1,
\end{equation}
where recall that $u_2\geq1$. Thus $x_1(T)<\kappa_1\leq1<\gamma$, and the final point $(\gamma,0)$ cannot be reached. Consequently, an extremal starting with a $Y$-segment can only end with a $Y$-segment, as we considered at the beginning.

If we consecutively apply (\ref{consecutive}) from the first switching point up to the last, we obtain
\begin{equation}
\label{Y_first_last}
\kappa_{2n}^2\kappa_{1}^2=\frac{(s+u_2)^{n-1}}{(s+u_1)^n}.
\end{equation}
Working as in the previous case we find
\begin{equation}
\label{Y_k1}
\kappa_{1,\pm}^{2}=\frac{c_{2}\pm\sqrt{c_{2}^{2}-4(s+u_{2})}}{2(s+u_{2})}=\frac{2}{c_{2}\mp\sqrt{c_{2}^{2}-4(s+u_{2})}}
\end{equation}
and
\begin{equation}
\label{Y_k_last}
\kappa_{2n,\pm}^{2}=\frac{c\pm\sqrt{c^{2}-4(s+u_{2})}}{2(s+u_{2})}=\frac{2}{c\mp\sqrt{c^{2}-4(s+u_{2})}},
\end{equation}
where $c_2=u_2+1$ and $c=u_2\gamma^2+1/\gamma^2$ as before. Since $c>c_{2}\Leftrightarrow(\gamma^{2}-1)(u_{2}\gamma^{2}-1)>0$ which is true, only the choice $\kappa_{2n,+}^{2}$ leads to a valid transcendental equation. Note that in order to exclude the choice $\kappa_{2n,-}^{2}$, one has to use in (\ref{Y_first_last}) the corresponding expression from the second equality in (\ref{Y_k_last}) and for $\kappa_{1,\pm}^{2}$ the expression from the first equality in (\ref{Y_k1}). If we use in (\ref{Y_first_last}) for $\kappa_{2n,+}^{2}$ the first equality in (\ref{Y_k_last}) and for $\kappa_{1,\pm}^{2}$ the second equality in (\ref{Y_k1}), we end up with the valid transcendental equation (\ref{transcendentalY}) for the ratio $s$, in the interval $0<s\leq (u_{2}-1)^{2}/4$, since $u_{2}-1<u_{2}\gamma^2-1/\gamma^2$ for $\gamma>1$. Having found $s$, the interswitching times $T_X, T_Y$ are given by (\ref{switch_X}), (\ref{switch_Y}), as above. The initial time interval $T_{I,2}^{\pm}$ can be calculated following the same procedure as before and the result is given in (\ref{time_in2}), while the final time interval $T_F$ is the same as in the previous case and is given in (\ref{time_fi}). The total duration $T_{2n}^{\pm}$ of
the trajectory with $2n$ switchings joining the points $(1,0)$ and $(\gamma,0)$ is
given by (\ref{time_even}), where $\pm$ corresponds to the choice of sign in (\ref{Y_k1}) for the first switching point.
\end{IEEEproof}

Using Theorem \ref{prop:time} we can find the times $T_{n}$ for a specific
target $(\gamma, 0)$ and compare them to obtain the minimum time. Some examples are given in the next section.

\begin{cor}
\label{uone}
For $u_2=1$, only extremals of the form $XY\ldots YXY$ with odd number of switchings are allowed.
\end{cor}

\begin{IEEEproof}
For $u=u_{2}=1$ the starting point $(1,0)$ is an equilibrium point of system (\ref{system1}), (\ref{system2}), so a trajectory cannot start with
a $Y$-segment. Note that for $u_2=1$, the upper bound for the allowed values of $s$ in the transcendental equation (\ref{transcendentalY}) is $(u_{2}-1)^{2}/4=0$, same as the lower bound.
\end{IEEEproof}

\begin{rem}
The major consequence of Theorem \ref{prop:time} is the possibility of odd-numbered extremals with more than one switchings for $u_2=1$. In the next section we present an example where such an extremal is actually the optimal solution. This kind of solution is not mentioned in any of the previous works \cite{Salamon09,Tsirlin11,Salamon12,Hoffmann13,Hoffmann15,Boldt16}.
\end{rem}

\begin{rem}
Solving numerically the transcendental equations for the ratio $s$ and then calculating the switching times using the formulas in Theorem \ref{prop:time}, is computationally more efficient than the numerical optimization over the switching times which is suggested in \cite{Boldt16}.
\end{rem}

\section{Examples}

\label{sec:examples}

\begin{table}[b]
\renewcommand{\arraystretch}{1.3}
\caption{Extremal Times}
\label{time_table}
\centering
\begin{tabular}{|c|c|c|c|c|}
\hline
& $\gamma=\sqrt{3}$  & $\gamma=\sqrt{3}$ & $\gamma=8$  & $\gamma=8$\\
& $u_2=1$ & $u_2=6.5$ & $u_2=1$ & $u_2=4$\\
\hline
$T_1^+$ & \bfseries 1.6784 & 1.4513 & 8.0159 & 7.9707\\
\hline
$T_3^+$ & - & - & \bfseries 7.3863 & 4.6189\\
\hline
$T_5^+$ & - & - & 9.5568 & -\\
\hline
$T_3^-$ & - & - & 9.7758 & 4.9845\\
\hline
$T_5^-$ & - & - & 9.5735 & -\\
\hline
$T_2^+$ & - & 1.8320 & - & 8.0452\\
\hline
$T_4^+$ & - & 2.5858 & - & 4.9982\\
\hline
$T_6^+$ & - & - & - & 5.7987\\
\hline
$T_8^+$ & - & - & - & 7.0651\\
\hline
$T_2^-$ & - & \bfseries 1.3888 & - & 4.8098\\
\hline
$T_4^-$ & - & 2.5387 & - & \bfseries 4.5458\\
\hline
$T_6^-$ & - & - & - & 5.6884\\
\hline
$T_8^-$ & - & - & - & 7.0496\\
\hline
\end{tabular}
\end{table}

In this section we illustrate the optimal solution described in Theorem \ref{prop:time} with several examples. For convenience we fix the lower control bound to $u_1=0.0002$ and consider four cases with the following \emph{realistic} values of $\gamma$ and $u_2$: (a) $\gamma=\sqrt{3}, u_2=1$, (b) $\gamma=\sqrt{3}, u_2=6.5$, (c) $\gamma=8, u_2=1$, (d) $\gamma=8, u_2=4$. For example, in the experiment \cite{Schaff11} where the fast but effectively adiabatic cooling of a trapped Bose-Einstein condensate is considered, the ratio of the initial to the final frequency is approximately $\omega_0/\omega_f=3$, corresponding to $\gamma=\sqrt{3}$. In Table \ref{time_table} we show for each of these cases the necessary times for the various extremals to reach the corresponding target point. These times result from Theorem \ref{prop:time}, by solving numerically the corresponding transcendental equations and subsequently using the formulas for the switching times. The absence of a solution is denoted by -, while the minimum time for each case is highlighted with bold.

\begin{figure}[t]
 \centering
		\begin{tabular}{cc}
     	\subfigure[$\ $$\gamma=\sqrt{3}, u_2=1$]{
	            \label{fig:ex1}
	            \includegraphics[width=.44\linewidth]{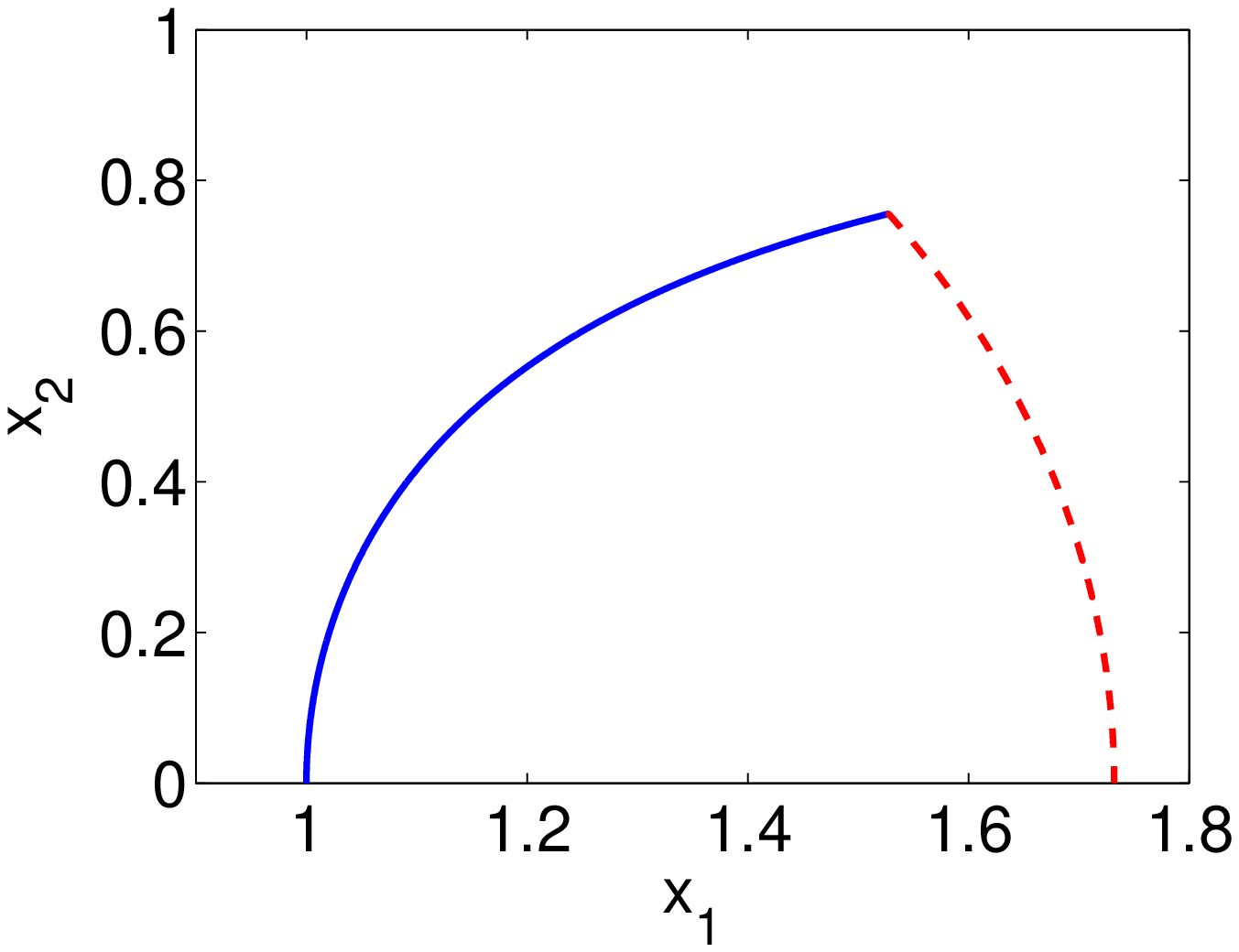}} &
	        \subfigure[$\ $$\gamma=\sqrt{3}, u_2=6.5$]{
	            \label{fig:ex2}
	            \includegraphics[width=.44\linewidth]{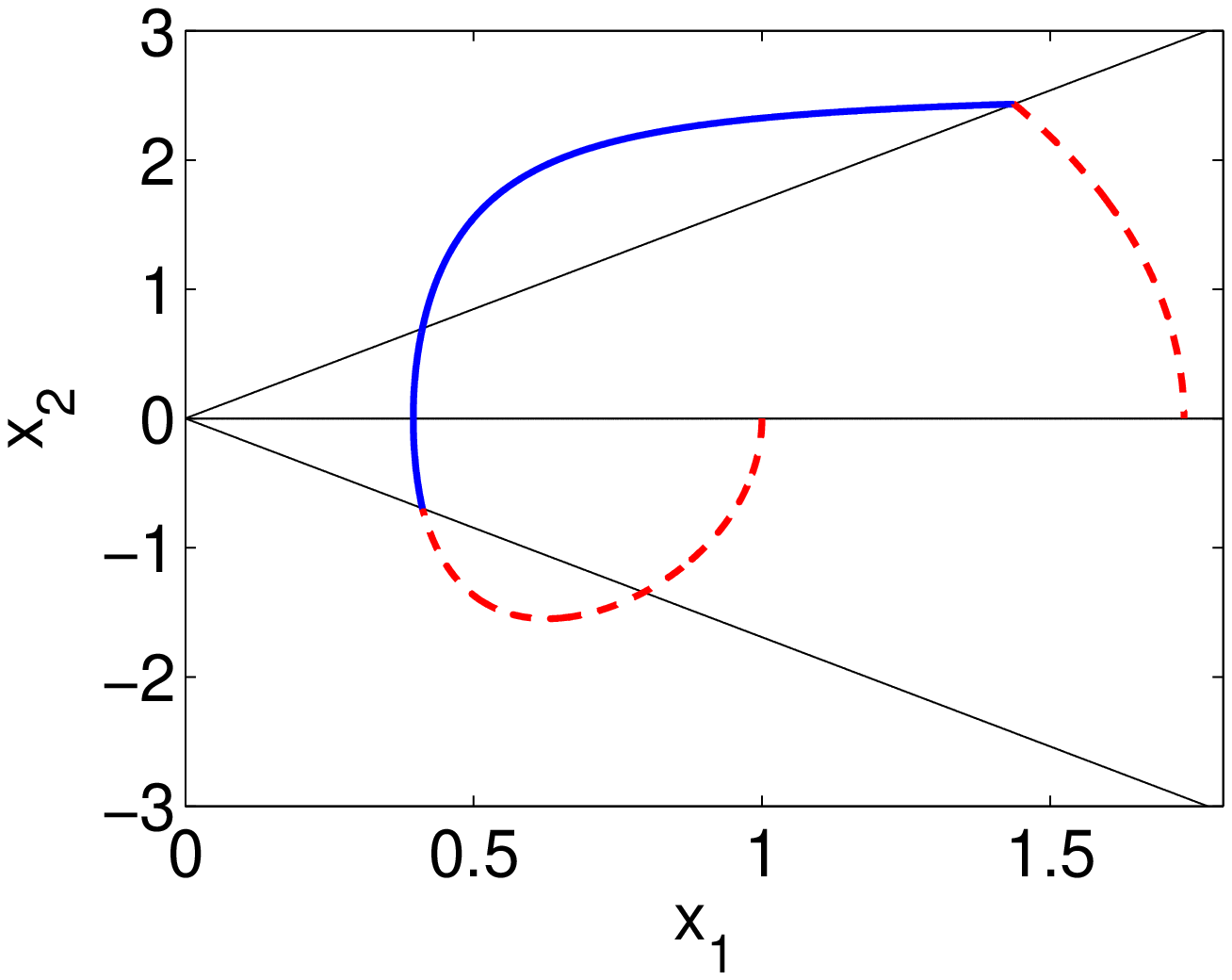}} \\
	        \subfigure[$\ $$\gamma=8, u_2=1$]{
	            \label{fig:ex3}
	            \includegraphics[width=.44\linewidth]{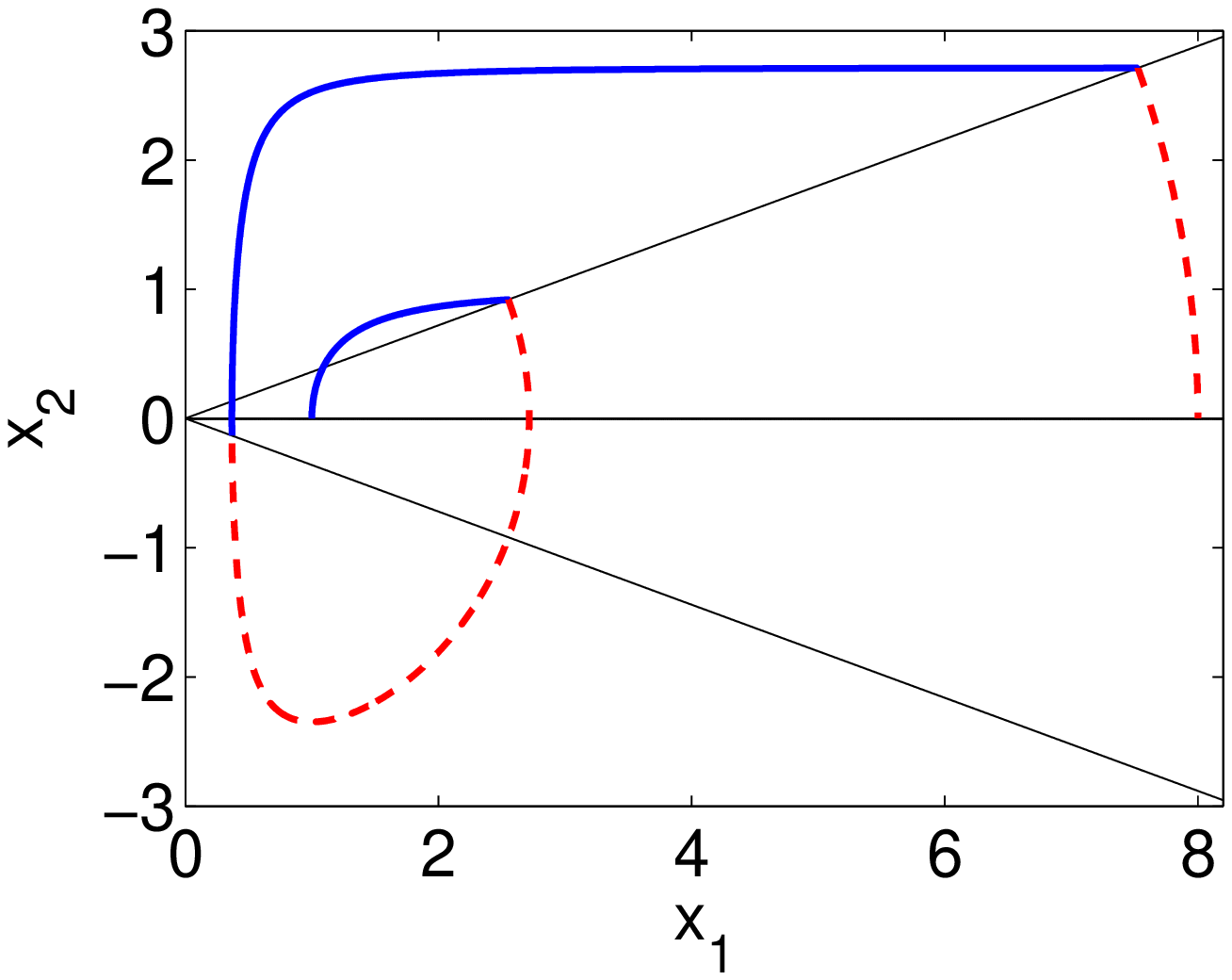}} &
			\subfigure[$\ $$\gamma=8, u_2=4$]{
	            \label{fig:ex4}
	            \includegraphics[width=.44\linewidth]{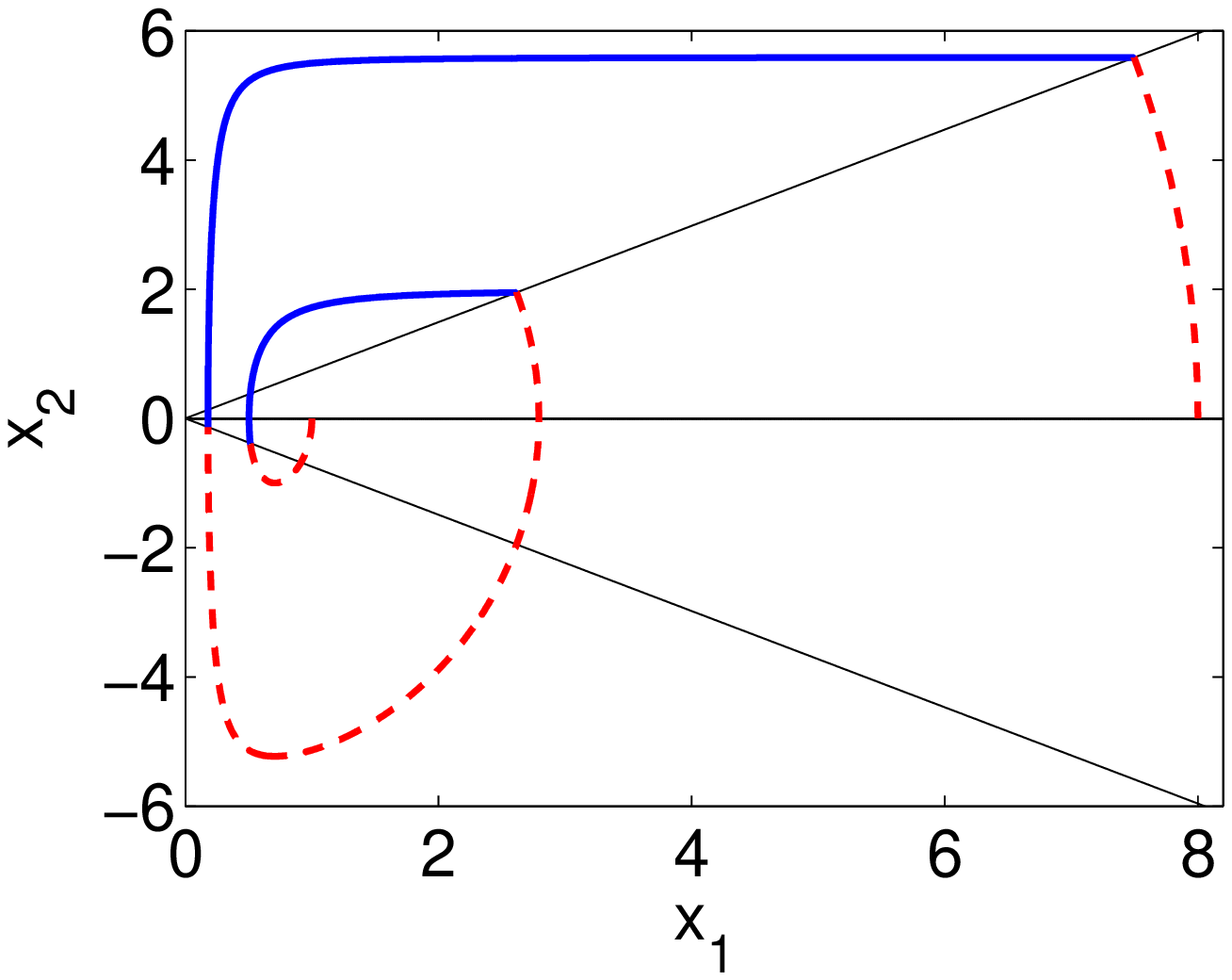}} \\
		\end{tabular}
\caption{Optimal trajectories for various values of the target $\gamma$ and the control upper bound $u_2$. The lower bound is fixed to $u_2=0.0002$ for all the cases.}
\label{fig:example}
\end{figure}

Observe that for the first case (first column), where $\gamma=\sqrt{3}, u_2=1$, there is only one extremal, and recall that the extremals starting with a $Y$-segment are excluded because $u_2=1$, as explained in Corollary \ref{uone}. The corresponding optimal trajectory $XY$ is depicted in Fig. \ref{fig:ex1}. For the next case (b), the target point is the same as before but the control upper bound has been increased to $u_2=6.5$. Table \ref{time_table} indicates that the minimum-time solution has two switchings, and the corresponding trajectory $YXY$ is shown in Fig. \ref{fig:ex2}. This kind of solution can be better understood if we adopt the point of view of Remark \ref{mech_analog} and interpret system (\ref{system1}), (\ref{system2}) as describing the one-dimensional Newtonian motion of a unit-mass particle, with $x_1, x_2$ corresponding to its position and velocity, respectively. If $u_2$ is large enough then the particle, instead of moving solely forward like in the $XY$ trajectory, can first approach $x_1=0$ sufficiently fast and then exploit the strong repulsive force $1/x_1^3$ to arrive faster at the target point. We identified this type of solution in \cite{Stefanatos11}, in the more general case where the control could also take negative values, but it was also verified numerically in \cite{Boldt16} for the restrictive case of positive controls.

We now move to the next example, where $\gamma=8, u_2=1$. As in case (a), the extremals starting with $Y$, having an even number of switching points, are excluded because of Corollary \ref{uone}. But, since the target point $\gamma$ is now larger than in the previous case, odd-numbered extremals with more than one switchings arise, which take advantage of the strong repulsive force close to $x_1=0$. In fact, the optimal solution is $XYXY$ with three switchings, as highlighted in Table \ref{time_table} and portrayed in Fig. \ref{fig:ex3}. We emphasize that this is a new kind of solution, absent from the previous works \cite{Salamon09,Tsirlin11,Salamon12,Hoffmann13,Hoffmann15,Boldt16}. These articles actually consider the case $0<\omega_f=\omega_1\leq\omega(t)\leq\omega_2=\omega_0$, which corresponds to $0<1/\gamma^4=u_1\leq u(t)\leq u_2=1$ in our terminology. For such restrictions in the frequency (stiffness) of the parametric oscillator, these papers conclude that the optimal solution has three ``jumps". This actually corresponds to our $XY$ solution with one (intermediate) switching, since in the number of jumps are included the changes at the initial and final times. In this language, our optimal solution with three (intermediate) switchings corresponds to a five-jump solution. Note that in our example we use $u_1=0.0002$, instead of $u_1=1/\gamma^4=2.44\cdot 10^{-4}$, but it can be verified that even in this case the optimal solution is again of the form $XYXY$ with three switchings. The crucial requirement for the validity of the comparison is $u_2=1$ on the upper bound, which assures that the even-numbered extremals are excluded, and not that on the lower bound \cite{Boldt16}. We finally mention that in our previous work \cite{Stefanatos11} we didn't identify odd-numbered solutions with more than one switchings since, for the control set that we considered there, $-u_1\leq u(t)\leq u_2$ with $u_1, u_2\geq 1$, this kind of extremals was excluded.

The last case that we examine has the same target point $\gamma=8$ but a larger control upper bound $u_2=4$. As we can observe from Table \ref{time_table} the optimal solution has four switchings, and the form YXYXY which is shown in Fig. \ref{fig:ex4}. Going back to the particle picture, we see that for these values of the parameters it is time-optimal to move back and forth twice in order to gain speed from the repulsive force, before reaching the target point. We close by pointing out that, using (\ref{zeta}) and (\ref{x}), the optimal trajectories shown in Fig. \ref{fig:example} can be easily displayed on the $z_1z_2$-plane, as in \cite{Boldt16}.

\section{Conclusion}

\label{sec:conclusion}

In this paper we used geometric optimal control to find the minimum necessary time for transitions between thermal states of the quantum parametric oscillator, and the corresponding optimal time-profile of oscillator's frequency. We considered the case where the frequency of the oscillator can take only real values, corresponding to nonnegative oscillator stiffness, and obtained a solution which has not appeared in any of the previous related works. The present work can find applications in several contexts, for example
to minimize the necessary time for the adiabatic stroke of a quantum heat engine and for the effectively adiabatic cooling of trapped atoms, reducing though the undesirable effect of random interactions with the environment, which is ubiquitous. In the future, we would like to apply a similar control theoretic approach to two systems that we have recently studied numerically: a quantum parametric oscillator with noise, which is a good model for a noisy quantum heat engine \cite{Stefanatos14PRE}, and two coupled oscillators with modulation in their coupling, which models pulsed cavity optomechanical cooling \cite{Stefanatos15}.



%

\begin{IEEEbiography}{Dionisis Stefanatos}
(M'11) was born in the Greek island of Cephalonia in 1977. He received the Diploma in Electrical Engineering with highest honors from NTU Athens and the PhD in Engineering Sciences from Harvard, where he was a co-recipient of the E. Jury best thesis award (2005). He has held postdoctoral positions at Harvard and Washington University in St. Louis, while currently is temporary lecturer in Physics at the Hellenic Army Academy in Athens. His research is focused on the study of control systems that arise from physical problems and especially quantum mechanical applications.
\end{IEEEbiography}







\end{document}